\documentclass[12pt]{article}

\makeatletter
\newfont{\footsc}{cmcsc10 at 8truept}
\newfont{\footbf}{cmbx10 at 8truept}
\newfont{\footrm}{cmr10 at 10truept}
\renewcommand{\ps@plain}{%
\renewcommand{\@oddfoot}{\footsc the electronic journal of combinatorics
  {\footbf 8} (2001), \#R00\hfil\footrm\thepage}}
\makeatother
\usepackage{amsmath,amsthm}
\usepackage{array}                  
\usepackage{pst-node}              % `PSTricks' v97 patch 14  <1999/12/23> (tvz)
\usepackage[enableskew, vcentermath]{youngtab}
\usepackage{url}                   % \url{email@address} sets email@address
\usepackage{caption2}               % This is for replacing Figure 1: with
\title{A \lq nice\rq\ bijection for a content formula for skew
  semistandard Young tableaux}

\author{Martin Rubey\\
\small Institut f\"ur Mathematik\\[-0.8ex]
\small Universit\"at Wien\\[-0.8ex]
\small %\url{http://www.univie.ac.at/\~{}rubey}\\[-0.8ex]
\url{http://www.univie.ac.at/~rubey}\\[-0.8ex]
\small \url{a9104910@unet.univie.ac.at}}

\date{\small Submitted: May 31, 2001;  Accepted: April 14, 2002.\\
\small MR Subject Classifications: 05E10}

\newtheorem{thm}{Theorem}[section]
\newtheorem{lem}[thm]{Lemma}

\theoremstyle{definition}
\newtheorem{dfn}[thm]{Definition}

\theoremstyle{remark}
\newtheorem*{rmk}{Remark}

\newcommand{\Dfn}[1]{\emph{#1}}                            % for definitions
\newcommand{\AlgR}{\textcircled{\footnotesize$\Rightarrow$}}
\newcommand{\AlgL}{\textcircled{\footnotesize$\Leftarrow$}}

\newcounter{seq}
\newcommand{\seq}[4][,]%             [separator] letter first last
{\setcounter{seq}{#3}%
\ifmmode%
#2{\arabic{seq}}#1\stepcounter{seq}#2{\arabic{seq}}#1\dots #1#2{#4}%
\else%
$#2{\arabic{seq}}$#1~\stepcounter{seq} $#2{\arabic{seq}}$#1 \dots #1~$#2{#4}$%
\fi}

\begin{document}
\maketitle

\begin{abstract}
  Based on Sch\"utzenberger's evacuation and a modification of jeu de
  taquin, we give a bijective proof of an identity connecting the
  generating function of reverse semistandard Young tableaux with
  bounded entries with the generating function of all semistandard
  Young tableaux. This solves Exercise 7.102 b of Richard Stanley's
  book `Enumerative Combinatorics~2'.
\end{abstract}

\section{Introduction}
The purpose of this article is to present a solution for Exercise
7.102 b of Richard Stanley's book `Enumerative
Combinatorics~2'~\cite{EC2}. There, Stanley asked for a `nice'
bijective proof of the identity
\begin{equation}\label{eq:rSSYT}
  \sum_{\substack{\text{$R$ reverse SSYT}\\
                  \text{of shape $\lambda/\mu$}\\
                  \text{with $R_{ij} \leq a+\mu_i-i$}}}
      q^{n(R)}
     = \bigg(\sum_{\substack{\text{$P$ SSYT}\\
                             \text{of shape $\lambda/\mu$}}} q^{n(P)}\bigg)
 \cdot \prod_{\rho\in\lambda/\mu}(1-q^{a+c(\rho)}),
\end{equation}
where $a$ is an arbitrary integer such that $a+c(\rho)>0$ for all cells
$\rho\in\lambda/\mu$.\footnote{In fact, this is the corrected version of the
  identity originally given in~\cite{EC2}, to be found at
  \url{http://www-math.mit.edu/~rstan/ec}.
  %\url{http://www-math.mit.edu/\~{}rstan/ec}. 
  Stanley took it from~\cite{BilleyJockuschStanley}, Theorem 3.1,
  where the formula is stated incorrectly, too.} Here, and in the
sequel, we use notation defined below:
\begin{dfn}
  A \Dfn{partition} is a sequence $\lambda=(\seq{\lambda_}{1}{r})$
  with $\seq[\geq]{\lambda_}{1}{r} > 0$, for some $r$.
  
  The \Dfn{Ferrers diagram} of a partition $\lambda$ is an array of
  cells with $r$ left-justified rows and $\lambda_i$ cells in row $i$.
  Figure~\ref{f:Diagrams}.a shows the Ferrers diagram corresponding to
  $(4, 3, 3, 1)$. We label the cell in the $i$\textsuperscript{th} row
  and $j$\textsuperscript{th} column of the Ferrers diagram of
  $\lambda$ by the pair $(i,j)$. Also, we write $\rho\in\lambda$, if
  $\rho$ is a cell of $\lambda$.
  
  A partition $\mu=(\seq{\mu_}{1}{s})$ is \Dfn{contained} in a
  partition $\lambda=(\seq{\lambda_}{1}{r})$, if $s\leq r$ and
  $\mu_i\leq\lambda_i$ for $i\in\{\seq{}{1}{s}\}$.
  
  The \Dfn{skew diagram} $\lambda/\mu$ of partitions $\lambda$ and
  $\mu$, where $\mu$ is contained in $\lambda$, consists of the cells
  of the Ferrers diagram of $\lambda$ which are not cells of the
  Ferrers diagram of $\mu$. Figure~\ref{f:Diagrams}.b shows the skew
  diagram corresponding to $(4, 3, 3, 1)/(2,2,1)$. The \Dfn{content}
  $c(\rho)$ of a cell $\rho=(i,j)$ of $\lambda/\mu$ is $j-i$.
  
  Given partitions $\lambda$ and $\mu$, a \Dfn{tabloid of shape $\lambda/\mu$}
  is a filling $T$ of the cells of the skew diagram $\lambda/\mu$ with
  non-negative integers. $T_\rho$ denotes the entry of $T$ in cell $\rho$. The
  \Dfn{norm} $n(T)$ of a tabloid $T$ is simply the sum of all entries of $T$.
  The \Dfn{content weight} $w_c(T)$ of a tabloid $T$ is
  $\sum_{\rho\in\lambda/\mu} T_\rho\cdot\left(a+c(\rho)\right)$, where $a$ is a
  given integer such that $a+c(\rho)>0$ for all cells $\rho\in\lambda/\mu$.
  
  A \Dfn{semistandard Young tableau of shape $\lambda/\mu$}, short
  \Dfn{SSYT}, is a tabloid $P$ such that the entries are weakly
  increasing along rows and strictly increasing along columns.
  
  A \Dfn{reverse semistandard Young tableau of shape $\lambda/\mu$} is
  a tabloid $R$ such that the entries are weakly decreasing along rows
  and strictly decreasing along columns. In Figure~\ref{f:Diagrams}.c
  a reverse SSYT of shape $(4,3,3,1)/(2,2,1)$ is shown.
\end{dfn}

\begin{figure}
  \begin{center}
    \begin{tabular}{ccc}
      \yng(4,3,3,1)
      & \young(::\hfil\hfil,::\hfil,:\hfil\hfil,\hfil)
      & \young(::33,::2,:30,4)\\
      a. Ferrers diagram
      & b. skew Ferrers diagram
      & c. reverse SSYT
    \end{tabular}
 \end{center}
  \caption{\label{f:Diagrams}}
\end{figure}

\section{A Bijective proof of Identity~\ref{eq:rSSYT}}
In fact, we will give a bijective proof of the following rewriting of
Identity~\ref{eq:rSSYT}:
\begin{align*}
  \sum_{\substack{\text{$P$ SSYT}\\
      \text{of shape $\lambda/\mu$}}} q^{n(P)}
  &= \bigg(\sum_{\substack{\text{$R$ reverse SSYT}\\
      \text{of shape $\lambda/\mu$}\\
      \text{with $R_{ij} \leq a+\mu_i-i$}}} q^{n(R)}\bigg)
  \cdot \prod_{\rho\in\lambda/\mu}\frac{1}{1-q^{a+c(\rho)}}\\
  &= \sum_{\substack{(R, T)\\
      \text{$R$ reverse SSYT}\\
      \text{of shape $\lambda/\mu$}\\
      \text{with $R_{ij} \leq a+\mu_i-i$,}\\
      \text{$T$ tabloid}\\
      \text{of shape $\lambda/\mu$}}} q^{n(R)} q^{w_c(T)}.
\end{align*}

So all we have to do is to set up a bijection that maps SSYT'x $P$
onto pairs $(R,T)$, where $R$ is a reverse SSYT with $R_{ij}\leq
a+\mu_i-i$ and $T$ is an arbitrary tabloid, such that
$n(P)=n(R)+w_c(T)$.

The bijection consists of two parts. The first step is a modification
of a mapping known as `evacuation', which consists of a special
sequence of so called `jeu de taquin slides'. An in depth description
of these procedures can be found, for example, in Bruce Sagan's Book
`The symmetric group'~\cite{SymmetricGroup}, Sections~3.9 and 3.11. We
use evacuation to bijectively transform the given SSYT $P$ in a
reverse SSYT $Q$ which has the same shape and the same norm as the
original one.

The second step of our bijection also consists of a sequence of --
modified -- jeu de taquin slides and bijectively maps a reverse SSYT
$Q$ onto a pair $(R,T)$ as described above. This procedure is very
similar to bijections discovered by Christian Krattenthaler, proving
Stanley's hook-content
formula.~\cite{Krattenthaler1998a,Krattenthaler1999a}

\begin{figure}[h]
  \begin{center}
    \begin{tabular}{ccccc}
      $\young(::01,::17,:149,299)$
      &$\Leftrightarrow$
      &$\young(::99,::74,:911,210)$
      &$\Leftrightarrow$
      &$\left(\,\young(::43,::22,:410,210)\,,
        \young(::00,::00,:102,001)\,\right)$\\\\
      $n(.)=43$
      & 
      &$n(.)=43$
      &
      &$n(.)=19, w_c(.)=24$
    \end{tabular}
  \end{center}
  \caption{}\label{f:Bijection}
\end{figure}

A complete example for the bijection can be found in the appendix.
There we chose $a=6$ and map the SSYT $P$ of shape $(4,4,4,3)/(2,2,1)$
on the left of Figure~\ref{f:Bijection} to the reverse SSYT $Q$ in the
middle of Figure~\ref{f:Bijection}, which in turn is mapped to the
pair on the right of Figure~\ref{f:Bijection}, consisting of a reverse
SSYT $R$, where the entry of the cell $\rho=(i,j)$ is less or equal to
$a+\mu_i-i$, and a tabloid $T$ so that $n(Q)=n(R)+w_c(T)$.

In the algorithm described below we will produce a filling of a skew
diagram step by step, starting with the `empty tableau' of the given
shape.
\begin{thm}\label{t:Evacuation}
  The following two maps define a correspondence between SSYT'x and reverse
  SSYT'x of the same shape $\lambda/\mu$ and the same norm:
\end{thm}
\begin{itemize}
\item[\AlgR] Given a SSYT $P$ of shape $\lambda/\mu$, produce a reverse SSYT
  $Q$ of the same shape and the same norm as follows:
\item[] Let $Q$ be the empty tableau of shape $\lambda/\mu$.
\item[] WHILE there is a cell of $P$ which contains an entry
  \begin{itemize}
  \item[] Let $e$ be the minimum of all entries of $P$. Among all
    cells $\tau$ with $P_\tau=e$, let $\rho=(i,j)$ be the cell which
    is situated most right.
  \item[] WHILE $\rho$ has a bottom or right neighbour in $P$ that
    contains an entry
    \begin{itemize}
    \item[] Denote the entry to the right of $\rho$ by $x$ and the
      entry below $\rho$ by $y$. We allow also that there is only an
      entry to the right or below $\rho$ and the other cell is missing
      or empty.
    \item[] If $x<y$, or there is no entry below $\rho$, then replace
      \begin{align*}
        \young(ex,y) \quad\text{ by }\quad \young(xe,y)\,,
      \end{align*}
      and let $\rho$ be the cell $(i,j+1)$.
    \item[] Otherwise, if $x\geq y$, or there is no empty to the right,
      replace
      \begin{align*}
        \young(ex,y) \quad\text{ by }\quad \young(yx,e)\,,
      \end{align*}
      and let $\rho$ be the cell $(i+1,j)$.
    \end{itemize}
  \item[] END WHILE.
  \item[] Put $Q_\rho$ equal to $e$ and delete the entry of the cell
    $\rho$ from $P$. Note that cells of $P$ which contain an entry
    still form a SSYT. In the proof below, $\rho$ will be called the
    cell where the jeu de taquin slide stops.
  \end{itemize}
\item[] END WHILE.
\item[\AlgL] Given a reverse SSYT $Q$ of shape $\lambda/\mu$, produce a SSYT
  $P$ of the same shape and the same norm as follows:
\item[] Let $P$ be the empty tableau of shape $\lambda/\mu$.
\item[] WHILE there is a cell of $Q$ which contains an entry
  \begin{itemize}
  \item[] Let $e$ be the maximum of all entries of $Q$. Among all
    cells $\tau$ with $Q_\tau=e$, let $\rho=(i,j)$ be the cell which
    is situated most left.
  \item[] Set $P_\rho=e$ and delete the entry of the cell $\rho$ from
    $Q$.
  \item[] WHILE $\rho$ has a top or left neighbour in $P$ that
    contains an entry
    \begin{itemize}
    \item[] Denote the entry to the left of $\rho$ by $x$ and the
      entry above $\rho$ by $y$. We allow also that there is only an
      entry to the left or above $\rho$ and the other cell is
      missing or empty.
    \item[] If $x>y$, or there is no entry above $\rho$, then replace
      \begin{align*}
        \young(:y,xe) \quad\text{ by }\quad \young(:y,ex)\,,
      \end{align*}
      and let $\rho$ be the cell $(i,j-1)$.
    \item[] Otherwise, if $x\leq y$, or there is no entry to the left,
      replace
      \begin{align*}
        \young(:y,xe) \quad\text{ by }\quad \young(:e,xy)\,,
      \end{align*}
      and let $\rho$ be the cell $(i-1,j)$.
    \end{itemize}
  \item[] END WHILE.
  \item[] The cells of $P$ which contain an entry now form a SSYT. In
    the proof below, $\rho$ will be called the cell where the jeu de
    taquin slide stops.
  \end{itemize}
\item[] END WHILE.
\end{itemize}
\begin{proof}
  Note that what happens during the execution of the inner loop of
  \AlgR\ (\AlgL) is a \Dfn{jeu de taquin forward (backward) slide}
  performed on $Q$ into the cell $\rho$, see Section~3.9
  of~\cite{SymmetricGroup}.
  
  First we have to show that \AlgR\ is well defined. I.e., we have to
  check that after each jeu de taquin forward slide, after the entry
  $e$ in the cell $\rho$ is deleted from $P$, the cells of $P$ which
  contain an entry form a SSYT as stated in the algorithm. This
  follows, because after either type of replacement in the inner loop
  the only possible violations of increase along rows and strict
  increase along columns in $P$ can only involve $e$ and the entries
  to its right and below. When the jeu de taquin forward slide is
  finished, $\rho$ is a bottom-right corner of $P$, hence after
  deleting the entry in $\rho$ no violations of increase or strict
  increase can occur.
  
  Next we show that \AlgR\ indeed produces a reverse SSYT. In fact, we
  even show that the tabloid defined by the cells of $Q$ which have
  been filled already, is a reverse SSYT at every stage of the
  algorithm.
  
  Clearly, every cell of $Q$ is filled with an entry exactly once.
  Furthermore, at the time the cell $\rho$ is filled, the cells in $Q$
  to the right and to the bottom of $\rho$ -- if they exist -- are
  filled already, otherwise $\rho$ would not be a bottom-right corner
  of $P$. Because the sequence of entries chosen is monotonically
  increasing, rows and columns of $Q$ are decreasing.
  
  So it remains to show that the columns of $Q$ are in fact strictly
  decreasing. Suppose that $\rho_1$ and $\rho_2$ are cells both
  containing the same minimal entry $e$, and $\rho_1$ is right of
  $\rho_2$.

\begin{figure}
  \begin{center}
    \Yautoscale1
    \begin{pspicture}(-0.5,-0.5)(4.5,4)\SpecialCoor
      \rput(0, 2.25){\rnode{r1}{\yng(1)}} \uput{9pt}[0](r1){$\rho_2$}
      \rput(1, 3.5) {\rnode{r2}{\yng(1)}} \uput{9pt}[0](r2){$\rho_1$}
      \rput(4, 0) {\rnode{r3}{\yng(1)}}
      \uput{9pt}[180](r3){$\rho_1^\prime$} \rput(4, 1)
      {\rnode{r4}{\yng(1)}} \uput{9pt}[90](r4){$\rho_2^\prime$} \rput(2,
      1.5) {\Rnode[href=0.5,vref=16pt]{c}{\yng(2,2)}}
      \nccurve[angleA=315, angleB=180]{r1}{c} \nccurve[angleA=315,
      angleB=90]{r2}{c} \nccurve[angleA=270, angleB=135]{c}{r3}
      \nccurve[angleA=0, angleB=135]{c}{r4}
   \end{pspicture}
  \end{center}
  \caption{}\label{f:PathsEvacuation}
\end{figure}
  When the jeu de taquin forward slide in \AlgR\ is performed into the
  cell $\rho_1$, the entry $e$ describes a path from $\rho_1$ to the
  cell where the slide stops, which we will denote by $\rho_1^\prime$.
  Similarly, we have a path from $\rho_2$ to a cell $\rho_2^\prime$.
  
  Now suppose $\rho_1^\prime$ is in the same column as, but below
  $\rho_2^\prime$, as depicted in Figure~\ref{f:PathsEvacuation}.
  Clearly, in this case the two paths would have to cross and we had
  the following situation:
  
  First, (the star is a placeholder for an entry we do not know)
  \begin{align*}
    \young(*c,zy) \quad\text{ would be replaced by }\quad
    \young(*y,zc).
  \end{align*}
  In this situation, $z$ would have to be smaller then $y$.

  Then, when the jeu de taquin forward slide into the cell $\rho_2$ is
  performed, the following situation would arise at the same four
  cells:
  \begin{align*}
    \young(cy,z*) \quad\text{ would have to be replaced by }\quad
    \young(yc,z*).
  \end{align*}
  But this cannot happen, because then $y$ would have to be strictly
  smaller than $z$.
  
  It can be shown in a very similar manner that \AlgL\ indeed produces
  a SSYT. We leave the details to the reader.
   
  Finally, we want to prove that \AlgL\ is inverse to \AlgR. Suppose
  that in \AlgR, a jeu de taquin forward slide into the cell $\rho$
  containing the entry $e$ is performed on $P$. Suppose that the slide
  stopped in $\rho^\prime$, $Q_{\rho^\prime}$ is set to $e$ and the
  entry in $\rho^\prime$ is deleted from $P$. Among the entries of
  $Q$, $e$ is maximal, because smallest entries are chosen first in
  \AlgR.  Furthermore, among those cells of $Q$ containing the entry
  $e$, the cell $\rho^\prime$ is most left. This follows, because the
  tabloid defined by the cells of $Q$ which have been filled already,
  is a reverse SSYT, and the paths defined by the jeu de taquin slides
  cannot cross, as we have shown above.
  
  It is straightforward to check that in this situation the jeu de
  taquin backward slide into $\rho^\prime$ performed on $P$ in \AlgL\
  stops in the original cell $\rho$. By induction we find that \AlgL\
  is inverse to \AlgR.
\end{proof}
The second step of the bijection is just as easy:
\begin{thm}\label{t:Bijection}
  The following two maps define a correspondence between reverse
  SSYT'x $Q$ to pairs $(R,T)$, where $R$ is a reverse SSYT with
  $R_{ij}\leq a+\mu_i-i$ and $T$ is an arbitrary tabloid, so that
  $n(Q) = n(R) + w_c(T)$, $Q$, $R$ and $T$ being of shape
  $\lambda/\mu$:
\end{thm}
\newcommand{\xplusone}{x+1} \newcommand{\xminusone}{x-1}
\newcommand{\yplusone}{y+1} \newcommand{\yminusone}{y-1}
\newcommand{\eminusone}{e-1} \Yboxdim25pt
\begin{itemize}
\item[\AlgR] Given a reverse SSYT $Q$ of shape $\lambda/\mu$, produce
  a pair $(R, T)$ as described above as follows:
\item[] Set $R=Q$ and set all entries of $T$ equal to $0$.
\item[] WHILE there is a cell $\tau=(i,j)$ such that
  $R_\tau>a+\mu_i-i$
  \begin{itemize}
  \item[] Let $e$ be maximal so that there is a cell $\tau$ with
    $R_\tau-\big(a+c(\tau)\big)=e$. Among all cells $\tau$ with
    $R_\tau-\big(a+c(\tau)\big)=e$, let $\rho=(i,j)$ be the cell which
    is situated most bottom. Set $R_\rho=e$.
  \item[] WHILE $e<R_{(i,j+1)}$ or $e\leq R_{(i+1,j)}$
    \begin{itemize}
    \item[] Denote the entry to the right of $\rho$ by $x$ and the
      entry below $\rho$ by $y$. We allow also that there is only a
      cell to the right or below $\rho$ and the other cell is missing.
    \item[] If $x-1>y$, or there is no cell below $\rho$, then replace
      \begin{align*}
        \young(ex,y) \quad\text{ by }\quad \young(\xminusone e,y)\,,
      \end{align*}
      and let $\rho$ be the cell $(i,j+1)$.
    \item[] Otherwise, if $y+1\geq x$, or there is no cell to the
      right, replace
      \begin{align*}
        \young(ex,y) \quad\text{ by }\quad \young(\yplusone x,e)\,,
      \end{align*}
      and let $\rho$ be the cell $(i+1,j)$.
    \end{itemize}
  \item[] END WHILE.
  \item[] Increase $T_\rho$ by one.
  \end{itemize}
\item[] END WHILE.
\item[\AlgL] Given a pair $(R, T)$ as described above, produce a
  reverse SSYT $Q$ of shape $\lambda/\mu$ as follows:
\item[] Set $Q=R$.
\item[] WHILE there is a cell $\tau=(i,j)$ such that $T_\tau \neq 0$
  \begin{itemize}
  \item[] Let $e$ be minimal so that there is a cell $\tau$ with
    $Q_\tau=e$ and $T_\tau\neq 0$. Among these cells $\tau$ let
    $\rho=(i,j)$ be the cell which is situated most right. Decrease
    $T_\rho$ by one.
  \item[] WHILE $e+a+c(\rho)>Q_{(i,j-1)}$ or $e+a+c(\rho)\ge Q_{(i-1,j)}$
                %$\rho$ is not a top-left corner
    \begin{itemize}
    \item[] Denote the entry to the left of $\rho$ by $x$ and the
      entry above $\rho$ by $y$. We allow also that there is only a
      cell to the left or above $\rho$ and the other cell is missing.
    \item[] If $y>x+1$, or there is no cell above $\rho$, then replace
      \begin{align*}
        \young(:y,xe) \quad\text{ by }\quad \young(:y,e\xplusone)\,,
      \end{align*}
      and let $\rho$ be the cell $(i,j-1)$.
    \item[] Otherwise, if $x\geq y-1$, or there is no cell to the
      left, replace
      \begin{align*}
        \young(:y,xe) \quad\text{ by }\quad \young(:e,x\yminusone)\,,
      \end{align*}
      and let $\rho$ be the cell $(i-1,j)$.
    \end{itemize}
  \item[] END WHILE.
  \item[] Increase $Q_\rho$ by $a+c(\rho)$.
  \end{itemize}
\item[] END WHILE.
\end{itemize}
\begin{rmk}
  Because of the obvious similarity to jeu de taquin slides, we will
  call what happens in the inner loop of \AlgR\ (\AlgL) a \Dfn{modified
    jeu de taquin (backward) slide} into $\rho$ performed on $R$
  ($Q$).
\end{rmk}
\begin{lem}\label{l:Paths}
  The two maps \ref{t:Bijection}.\AlgR\ and \ref{t:Bijection}.\AlgL\ are
  well defined. I.e. the tabloid $R$ produced by \AlgR\ is indeed a
  reverse SSYT with $R_{ij}\leq a+\mu_i-i$ and the tabloid $Q$
  produced by \AlgL\ is indeed a reverse SSYT. Also, the equation $n(Q)
  = n(R) + w_c(T)$ holds.
  
  Furthermore, the following statement is true: Suppose that \AlgR\
  performs a modified jeu de taquin slide on $R$ into a cell $\rho_1$
  with $R_{\rho_1}=e$.  After this, suppose that another modified jeu
  de taquin slide on $R$ into a cell $\rho_2$ with the same entry $e$
  is performed.  Let $\rho_1^\prime$ and $\rho_2^\prime$ be the cells
  where the slides stop. Then $\rho_1^\prime$ is left of
  $\rho_2^\prime$ or $\rho_1^\prime=\rho_2^\prime$. A corresponding
  statement holds for Algorithm~\ref{t:Bijection}.\AlgL.
\end{lem}
\begin{proof} 
  First of all, we have to prove that Algorithm~\ref{t:Bijection}.\AlgR\
  terminates. We required that $a+c(\tau)>0$ for all cells $\tau$, which
  implies that every time when we replace the entry in cell $\rho$ by $e$ (see
  the beginning of the outer loop of the algorithm) we decrease
  $\max_{\tau=(i,j)} (R_\tau-a-\mu_i+i)$. It is easy to see that this maximum
  is never increased in the subsequent steps of the algorithm.

  It is easy to check that after every type of replacement
  within the modified jeu de taquin slides, the validity of the equation $n(Q)
  = n(R) + w_c(T)$ is preserved.
  
  So it remains to show that after every modified jeu de taquin slide
  of \AlgR, the resulting filling $R$ of $\lambda/\mu$ is in fact a
  reverse SSYT: We have that $Q_\tau-\big(a+c(\tau)\big)=e$ is maximal
  at the very left of $\lambda/\mu$, because rows are decreasing in
  $Q$. Therefore, when $Q_\tau>a+\mu_i-i$, as required for the
  execution of the outer loop of \AlgR, we have
  \begin{align*}
    e=Q_\tau-\big(a+c(\tau)\big)>a+\mu_i-i-(a+\mu_i+1-i)=-1,
  \end{align*}
  so $e$ is non-negative.  Furthermore, after either type of
  replacement during the modified jeu de taquin slide, the only
  possible violations of decrease along rows or strict decrease along
  columns can involve only the entry $e$ and the entries to the right
  and below.  By induction, $R$ must be a reverse SSYT.
  
  The second statement of the lemma is shown with an argument similar
  to that used in the proof of Theorem~\ref{t:Evacuation}.

\begin{figure}
  \begin{center}
    \Yautoscale1
    \begin{pspicture}(-0.5,-0.5)(4.5,4)\SpecialCoor
      \rput(0, 2) {\rnode{r1}{\yng(1)}} \uput{9pt}[0](r1){$\rho_1$}
      \rput(1, 3.5) {\rnode{r2}{\yng(1)}} \uput{9pt}[0](r2){$\rho_2$}
      \rput(2.5, 0) {\rnode{r3}{\yng(1)}}
      \uput{9pt}[90](r3){$\rho_2^\prime$} \rput(4, 1)
      {\rnode{r4}{\yng(1)}} \uput{9pt}[90](r4){$\rho_1^\prime$} \rput(2,
      1.75){\Rnode[href=-0.5,vref=5pt]{c}{\yng(2,2)}}
      \nccurve[angleA=315, angleB=180]{r1}{c} \nccurve[angleA=315,
      angleB=90]{r2}{c} \nccurve[angleA=270, angleB=135]{c}{r3}
      \nccurve[angleA=0, angleB=135]{c}{r4}
    \end{pspicture}
  \end{center}
  \caption{}\label{f:Paths}
\end{figure}
  When the jeu de taquin forward slide in \AlgR\ is performed into the
  cell $\rho_1$, the entry $e$ describes a path from $\rho_1$ to the
  cell $\rho_1^\prime$, where the slide stops.  Similarly, we have a
  path from $\rho_2$ to $\rho_2^\prime$. We conclude that, if
  $\rho_1^\prime$ were strictly to the right of $\rho_2^\prime$, that
  these paths would have to cross. (See Figure~\ref{f:Paths}).
  Hence we had the following situation:
  
  First, (the star is a placeholder for an entry we do not know)
  \begin{align*}
    \young(*z,ex) \quad\text{ would be replaced by }\quad
    \young(*z,\xminusone e).
  \end{align*}
  In this situation, $x$ would have to be strictly smaller then $z$.
  
  Then, when the modified jeu de taquin slide into $\rho_2$ is
  performed, the following situation would arise at the same four
  cells:
  \begin{align*}
    \young(ez,\xminusone *)\quad\text{ would have to be replaced by
      }\quad 
    \young(xz,e*).
  \end{align*}
  But this cannot happen, because then $x$ would have to be at least
  as big as $z$ is.
  
  The corresponding statement for Algorithm~\ref{t:Bijection}.\AlgL\ is
  shown similarly.
\end{proof}

\begin{proof}[Proof of Theorem~\ref{t:Bijection}]
  It remains to show, that \AlgR\ and \AlgL\ are inverse to each other.
  This is pretty obvious considering the lemma:
  
  Suppose that the pair $(R,T)$ is an intermediate result obtained after a
  modified jeu de taquin slide into the cell $\rho$.  After this,
  $T_{\rho^\prime}$ is increased, where $\rho^\prime$ is the cell where the
  slide stopped. Then the entry in $\rho^\prime$ must be among the smallest
  entries of $R$, so that $T_{\rho^\prime}\neq 0$, because the sequence of
  $e$'s in the cells chosen for the modified jeu de taquin slides is
  monotonically decreasing.  If there is more than one cell $\rho$ which
  contains a minimal entry of $R$ and satisfies $T_\rho\neq 0$, the lemma
  asserts that the right-most cell was the last cell chosen for the modified
  jeu de taquin slide \AlgR.
  
  Hence it is certain that the right-most cell containing a minimal
  entry as selected before the modified jeu de taquin slide of \AlgL\
  is $\rho^\prime$.  
  It is easy to check, that the replacements done in \AlgL\ are exactly
  inverse to those in \AlgR. For example, suppose the following
  replacement is performed in \AlgR:
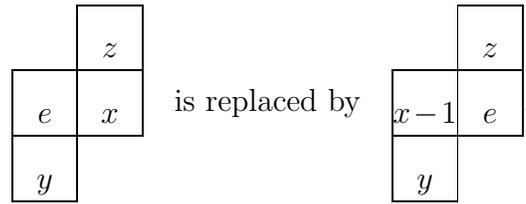
  \begin{equation*}
    \young(:z,ex,y)\quad\text{is replaced by}\quad
    \young(:z,\xminusone e,y).
  \end{equation*}
  Then we had $x-1>y$ and, because of strictly decreasing columns,
  $z>x$.  Therefore, in \AlgL, this is reversed and we end up with the
  original situation.
  
  Similarly, we can show that \AlgR\ is inverse to \AlgL, too.
\end{proof}
\appendix
\renewcommand\thesection{Appendix \Alph{section}:}
\section{Step by step example}
This appendix contains a complete example for the algorithms described above
for a SSYT of shape $(4,4,4,3)/(2,2,1)$ and $a=6$.

First the SSYT $P$ on the left of Figure~\ref{f:Bijection} is transformed into
the reverse SSYT $Q$ in the middle of Figure~\ref{f:Bijection} using
Algorithm~\ref{t:Evacuation}.\AlgR. The example has to be read in the following
way: Each pair $(P,Q)$ in the table depicts an intermediate result of the
algorithm. The cell of $P$ containing the encircled entry is the cell into
which the next jeu de taquin slide is performed. The jeu de taquin path is
indicated by the line in $Q$.
\Yautoscale1
\psset{framesep=0.5pt}
\newcommand{\ci}[1]{\pscirclebox[linewidth=0.2pt]{#1}}
\newcommand{\ciz}{\ci 0}\newcommand{\cio}{\ci 1}
\newcommand{\cit}{\ci 2}\newcommand{\cif}{\ci 4}
\newcommand{\cis}{\ci 7}\newcommand{\cie}{\ci 8}
\newcommand{\cin}{\ci 9}
\newcommand{\fb}[1]{\psframebox[linewidth=0.2pt]{#1}}
\newcommand{\fbz}{\fb 0}\newcommand{\fbo}{\fb 1}
\newcommand{\fbf}{\fb 4}\newcommand{\fbn}{\fb 9}
\newcommand{\cn}[1]{\circlenode[linewidth=0.2pt]{cn}{#1}}
\newcommand{\nn}{\Rnode[vref=4pt]{nn}{\hfil}}
\newcommand{\nnn}{\Rnode[vref=4pt]{nnn}{\hfil}}
\newcommand{\nnnn}{\Rnode[vref=4pt]{nnnn}{\hfil}}
\newcommand{\cnz}{\cn 0}\newcommand{\cno}{\cn 1}
\newcommand{\cnt}{\cn 2}\newcommand{\cnf}{\cn 4}
\newcommand{\cns}{\cn 7}\newcommand{\cnn}{\cn 9}
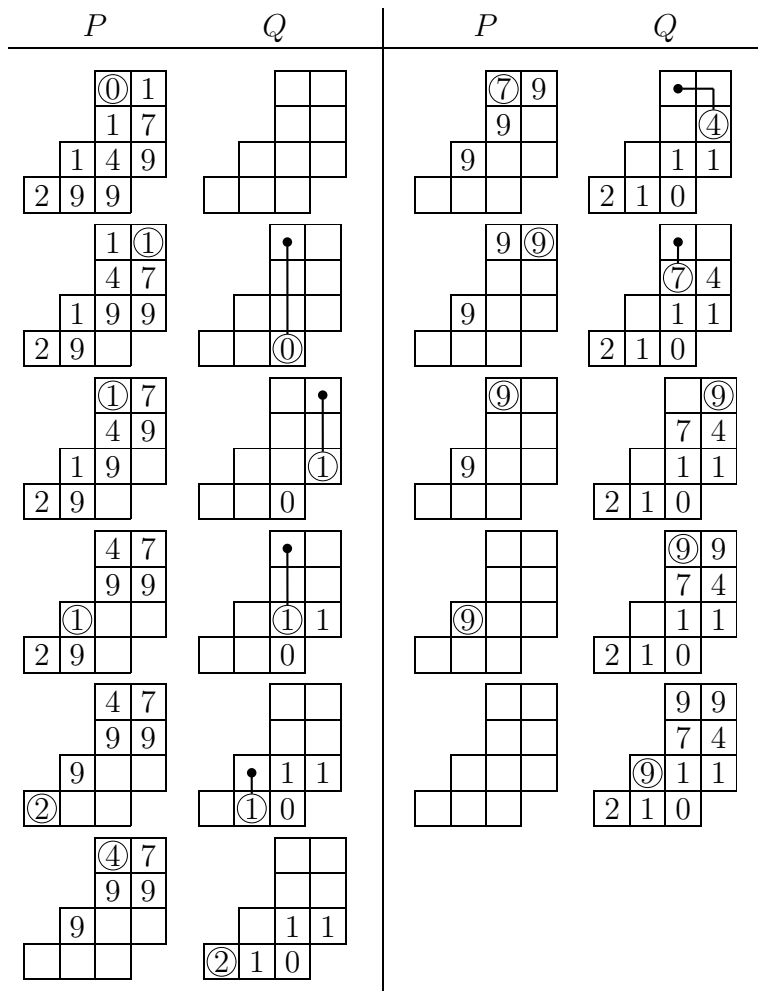
\begin{center}
\begin{tabular}[b]{cc@{\quad\vline\quad}cc}
  $P$ & $Q$ & $P$ & $Q$ \\
  \hline & & & \\[-6pt]
  \young(::\ciz 1,::17,:149,299) &
  \young(::\hfil\hfil,::\hfil\hfil,:\hfil\hfil\hfil,\hfil\hfil\hfil) &

  \young(::\cis 9,::9\hfil,:9\hfil\hfil,\hfil\hfil\hfil) &
  \young(::\nn\nnn,::\hfil\cnf,:\hfil 11,210)
  \ncline[arrows=*-]{nn}{nnn}
  \ncline{nnn}{cn}\\

  \young(::1\cio,::47,:199,29\hfil) &
  \young(::\nn\hfil,::\hfil\hfil,:\hfil\hfil\hfil,\hfil\hfil\cnz)
  \ncline[arrows=*-]{nn}{cn} &

  \young(::9\cin,::\hfil\hfil,:9\hfil\hfil,\hfil\hfil\hfil) &
  \young(::\nn\hfil,::\cns 4,:\hfil 11,210)
  \ncline[arrows=*-]{nn}{cn}\\

  \young(::\cio 7,::49,:19\hfil,29\hfil) &
  \young(::\hfil\nn,::\hfil\hfil,:\hfil\hfil\cno,\hfil\hfil 0) 
  \ncline[arrows=*-]{nn}{cn} &

  \young(::\cin\hfil,::\hfil\hfil,:9\hfil\hfil,\hfil\hfil\hfil) &
  \young(::\hfil\cin,::74,:\hfil 11,210)\\

  \young(::47,::99,:\cio\hfil\hfil,29\hfil) &
  \young(::\nn\hfil,::\hfil\hfil,:\hfil\cno 1,\hfil\hfil 0) 
  \ncline[arrows=*-]{nn}{cn} &

  \young(::\hfil\hfil,::\hfil\hfil,:\cin\hfil\hfil,\hfil\hfil\hfil) &
  \young(::\cin 9,::74,:\hfil 11,210)\\

  \young(::47,::99,:9\hfil\hfil,\cit\hfil\hfil) &
  \young(::\hfil\hfil,::\hfil\hfil,:\nn 11,\hfil\cno 0) 
  \ncline[arrows=*-]{nn}{cn} &

  \young(::\hfil\hfil,::\hfil\hfil,:\hfil\hfil\hfil,\hfil\hfil\hfil) &
  \young(::99,::74,:\cin 11,210)\\

  \young(::\cif 7,::99,:9\hfil\hfil,\hfil\hfil\hfil) &
  \young(::\hfil\hfil,::\hfil\hfil,:\hfil 11,\cit 10)
\end{tabular}
\end{center}
The inverse transformation \AlgL\ of the reverse SSYT $Q$ into the SSYT
$P$ can be traced in the same table, we only have to start at the
right bottom, where the tableau $P$ is empty, and work our way upwards
to the top left of the table. Note that the jeu de taquin paths are
the same.

In the second step of the bijection, this reverse SSYT $Q$ is mapped
onto a pair $(R,T)$, where $R$ is a reverse SSYT with $R_{ij}\leq
a+\mu_i-i$, $T$ is a tabloid and $n(Q)=n(R)+w_c(T)$. 

\begin{figure}[h]
  \begin{center}
    \begin{tabular}{cc}
      \young(7,6,4,2) &
      \young(::89,::78,:567,345) \\
      a. $a+\mu_i-i$ & b. The tabloid with \\
      & entries $a+c(\rho)$
    \end{tabular}
  \end{center}
  \caption{}\label{f:ContentTabloid}
\end{figure}

First, the algorithm initialises $R$ to $Q$ and sets all entries of
$T$ to zero. Using modified jeu de taquin slides, $R$ is then
transformed into a reverse SSYT where the entries are bounded as
required.  First the algorithm checks whether there are still cells in
$R$ which are too large. For reference, we give the relevant bounds in
Figure~\ref{f:ContentTabloid}.a. Then, for selecting the cell into
which the modified jeu de taquin slide is performed, we need to
calculate $R_\rho - \left(a + c(\rho)\right)$. Again, for reference we
display these values for each cell in Figure~\ref{f:ContentTabloid}.b.

Each row of the table below depicts an intermediate result of
Algorithm~\ref{t:Bijection}.\AlgR. The cells containing the encircled
entry are the cells into which the modified jeu de taquin slide will
be performed, the cells containing the boxed entry indicate, where the
last modified jeu de taquin slide stopped. In the third column the jeu
de taquin path for the selected cell is indicated.
\begin{center}
\begin{tabular}{cc}
  $R$ & $T$\\
  \hline\\[-6pt]
  \young(::99,::74,:\cin 11,210) & \young(::00,::00,:000,000) 
\\
  \young(::\cin 9,::74,:\fbf 11,210) & \young(::00,::00,:100,000) 
\\
  \young(::85,::\cis 2,:41\fbo,210) & \young(::00,::00,:101,000) 
\\
  \young(::\cie 5,::22,:411,21\fbz) & \young(::00,::00,:101,001) 
\\
  \young(::43,::22,:41\fbz,210) & \young(::00,::00,:102,001)
\end{tabular}
\begin{tabular}{c}
jeu de taquin path\\
\hline
\\[23pt]
  \young(::\hfil\hfil,::\hfil\hfil,:\nn\hfil\hfil,\hfil\hfil\hfil)
  \ncline[arrows=*-*]{nn}{nn}
\\
  \young(::\nn\nnn,::\hfil\hfil,:\hfil\hfil\nnnn,\hfil\hfil\hfil)
  \ncline[arrows=*-]{nn}{nnn}
  \ncline[arrows=-*]{nnn}{nnnn}
\\
  \young(::\hfil\hfil,::\nn\hfil,:\hfil\hfil\hfil,\hfil\hfil\nnn)
  \ncline[arrows=*-*]{nn}{nnn}
\\
  \young(::\nn\nnn,::\hfil\hfil,:\hfil\hfil\nnnn,\hfil\hfil\hfil)
  \ncline[arrows=*-]{nn}{nnn}
  \ncline[arrows=-*]{nnn}{nnnn}
\\[29pt]
\end{tabular}
\end{center}
Again, the inverse transformation \AlgL\ can be traced in the same
table, starting at the bottom, moving upwards. Now the cells
containing the boxed entry are the cells into which the next modified
jeu de taquin slide will be performed, the cells containing the
encircled entry indicate where the last slide stopped. Of course, the
jeu de taquin paths are the same as for \AlgR.
\section{A complete matchup for SSYT'x of\\ shape $(3,2)/(1)$ with norm $5$,
  where $a=2$} 

In the table below you find a complete matchup for SSYT'x of shape $(3,2)/(1)$
with norm $5$ where $a=2$. The first column contains all SSYT'x of shape
$(3,2)/(1)$ and norm $5$. In the second column, the corresponding reverse
SSYT'x obtained by evacuation are displayed. Finally, in columns three and
four, the results of Algorithm~\ref{t:Bijection}.\AlgR\ can be found.

This table was produced with a Common-LISP-implementation of the algorithms
above, which can be found on the author's
homepage.\footnote{%\url{http://www.mat.univie.ac.at/\~{}rubey/biject.lisp}
\url{http://www.mat.univie.ac.at/~rubey/biject.lisp}
}
\begin{center}
\setlength{\extrarowheight}{12pt}
\begin{equation*}
\begin{array}{c!{\Leftrightarrow}c!{\Leftrightarrow\bigg(}c!{,}c!{\bigg)}}
\multicolumn{1}{c}{P\phantom{\Leftrightarrow}}
&\multicolumn{1}{c}{Q\phantom{\Leftrightarrow\bigg(}}
&\multicolumn{1}{c}{R\phantom{,}}
&\multicolumn{1}{c}{T\phantom{\bigg)}}\\
\hline
\young(:00,05) & \young(:50,00) & \young(:20,00) & \young(:10,00) \\[12pt] \hline
\young(:00,14) & \young(:40,10) & \young(:10,00) & \young(:10,10) \\[12pt] \hline
\young(:00,23) & \young(:30,20) & \young(:10,00) & \young(:00,21) \\[12pt] \hline
\young(:01,04) & \young(:10,40) & \young(:10,00) & \young(:00,40) \\[12pt] \hline
\young(:01,13) & \young(:30,11) & \young(:20,00) & \young(:00,11) \\[12pt] \hline
\young(:01,22) & \young(:20,21) & \young(:20,00) & \young(:00,11) \\[12pt] \hline
\young(:02,03) & \young(:20,30) & \young(:20,00) & \young(:00,30) \\[12pt] \hline
\young(:02,12) & \young(:22,10) & \young(:22,00) & \young(:00,10) \\[12pt] \hline
\young(:11,03) & \young(:11,30) & \young(:11,00) & \young(:00,30) 
\end{array}
\end{equation*}
\begin{equation*}
\begin{array}{c!{\Leftrightarrow}c!{\Leftrightarrow\bigg(}c!{,}c!{\bigg)}}
\multicolumn{1}{c}{P\phantom{\Leftrightarrow}}
&\multicolumn{1}{c}{Q\phantom{\Leftrightarrow\bigg(}}
&\multicolumn{1}{c}{R\phantom{,}}
&\multicolumn{1}{c}{T\phantom{\bigg)}}\\
\hline
\young(:11,12) & \young(:21,11) & \young(:21,00) & \young(:00,01) \\[12pt] \hline
\young(:03,02) & \young(:32,00) & \young(:10,00) & \young(:01,00) \\[12pt] \hline
\young(:03,11) & \young(:31,10) & \young(:11,00) & \young(:00,11) \\[12pt] \hline
\young(:12,02) & \young(:21,20) & \young(:21,00) & \young(:00,20) \\[12pt] \hline
\young(:04,01) & \young(:41,00) & \young(:11,00) & \young(:10,00)
\end{array}
\end{equation*}
\end{center}
\newcommand{\cocoa} {\mbox{\rm C\kern-.13em o\kern-.07em C\kern-.13em
  o\kern-.15em A}}
\providecommand{\bysame}{\leavevmode\hbox to3em{\hrulefill}\thinspace}
\providecommand{\MR}{\relax\ifhmode\unskip\space\fi MR }
% \MRhref is called by the amsart/book/proc definition of \MR.
\providecommand{\MRhref}[2]{%
  \href{http://www.ams.org/mathscinet-getitem?mr=#1}{#2}
}
\providecommand{\href}[2]{#2}

\end{document}